\documentclass[aop,MSNbibl,citesort,dvips]{arximspdf}

\doi{10.1214/10-AOP632}
\volume{40}
\issue{2}
\pubyear{2012}
\firstpage{437}
\lastpage{458}

\makeatletter
\newtheorem{theo}{Theorem}[section]
\newtheorem{theorem}{Theorem}[section]
\newtheorem{cor}[theorem]{Corollary}
\newtheorem{prop}[theorem]{Proposition}
\newproclaim{rem}{Remark}[section]
\def\T{\Theta}

\def\C{{\mathbb C}}
\def\Z{{\mathbb Z}}
\def\R{{\mathbb R}}
\def\P{{\mathcal P}}

\def\Q{{\mathcal Q}}
\def\qhat{{\mathcal R}}
\def\L{{\mathcal L}}
\def\F{{\mathcal F}}
\def\T{{\mathcal T}}
\def\G{{\mathcal G}}
\def\A{{\mathcal A}}
\makeatother

\begin{document}
\begin{frontmatter}

\title{Directed polymers and the quantum Toda lattice}
\runtitle{Directed polymers and the quantum Toda lattice}

\begin{aug}
\author[A]{\fnms{Neil} \snm{O'Connell}\corref{}\ead[label=e1]{n.m.o-connell@warwick.ac.uk}}
\runauthor{N. O'Connell}
\affiliation{University of Warwick}
\address[A]{Mathematics Institute\\
University of Warwick\\
Coventry CV4 7AL\\
United Kingdom\\
\printead{e1}} 
\end{aug}

\received{\smonth{10} \syear{2009}}
\revised{\smonth{10} \syear{2010}}

%
\begin{abstract}
We characterize the law of the partition function
of a Brownian directed polymer model in terms of a diffusion process
associated with the quantum Toda lattice.
The proof is via a multidimensional generalization of a theorem of
Matsumoto and Yor
concerning exponential functionals of Brownian motion. It is based on a mapping
which can be regarded as a geometric variant of the RSK correspondence.
\end{abstract}

%

\end{frontmatter}

\section{Introduction}

Let $B_1(t),B_2(t),\ldots ,B_N(t), t\ge0$, be a collection of
independent standard
one-dimensional Brownian motions and write $B_i(s,t)=B_i(t)-B_i(s)$ for
$s\le t$.
Let $\beta\in\R$, $t\ge0$, and consider the random variable
\[
Z^N_t(\beta) = \int_{0<s_1<\cdots<s_{N-1}<t}
e^{\beta(B_1(s_1)+B_2(s_1,s_2)+\cdots+ B_N(s_{N-1},t))}\,ds_1\cdots \,ds_{N-1}.
\]
This is the partition function for a model of a $1+1$ dimensional
directed polymer
in a random environment
which has been introduced and studied in the papers \cite{om,oy2,sv}.
The free energy density is given explicitly by
\[
\lim_{N\to\infty} \frac1N\log Z^N_N(\beta)=
\inf_{t>0}[ \beta^2t - \Psi(t)]-\log\beta^2
\]
almost surely,
where $\Psi(z)=\Gamma'(z)/\Gamma(z)$.
The law of $Z^N_t(\beta)$ is well understood in the zero temperature
limit $\beta\to\infty$,
where it has close connections with random matrices. Define
%
\begin{eqnarray}\label{mdef}
M^N_t &=& \lim_{\beta\to\infty} \frac1\beta\log Z^N_t(\beta)
\nonumber
\\[-8pt]
\\[-8pt]
&=&
\max_{0\le s_1\le\cdots\le s_{N-1}\le t}
\bigl(B_1(s_1)+B_2(s_1,s_2)+\cdots+ B_N(s_{N-1},t)\bigr).
\nonumber
\end{eqnarray}
Note that, by Brownian scaling, the law of $M^N_t/\sqrt{t}$ is
independent of $t$.
%
\begin{theorem} The random variable $M^N_1$ has the same distribution
as the largest eigenvalue of a $N\times N$ GUE random matrix \cite{bar,gtw}.
In fact \cite{bj,oy}, the stochastic
process $(M^N_t, t\ge0)$ has the same law as the largest eigenvalue of
a standard Hermitian
Brownian motion, that is, it has the same law as the first coordinate
of a Brownian motion conditioned
(in the sense of Doob) never to exit the Weyl chamber $C_N=\{x\in\R
^N\dvtx  x_1>\cdots>x_N\}$,
started from the origin. This is a diffusion process in $\overline
{C}_N$ with infinitesimal
generator
$ \Delta/2 + \nabla\log h \cdot\nabla$
where
%
\begin{equation}\label{h}
h(x)=\prod_{1\le i<j\le N} (x_i-x_j).
\end{equation}
\end{theorem}

This connection with random matrices
yields very precise information concerning the distribution and
asymptotic behavior of $M^N$ when
$N$ is large. For example, it follows that
\[
\lim_{N\to\infty} P ( M^N_N \le2N + xN^{1/3}  ) = F_2(x),
\]
where $F_2$ is the Tracy--Widom distribution \cite{tw}.

In this paper we obtain an analogue of Theorem 1.1 for the stochastic process
$(\log Z^N_t(\beta), t>0)$. We will show that, for each $\beta>0$,
this process
has the same law as the first coordinate of a diffusion process in $\R
^N$ which is closely
related to the quantum Toda lattice. This yields an analytic description
of the law of $Z^N_t(\beta)$ which should provide a good starting
point for further
asymptotic analysis.

\section{The quantum Toda lattice}

The quantum Toda lattice is a quantum integrable system with Hamitonian
given by
the Schr\"{o}dinger operator
%
\begin{equation}
H=\sum_{i=1}^N\frac{\partial^2}{\partial x_i^2} - 2 \sum
_{i=1}^{N-1} e^{x_{i+1}-x_i} .
\end{equation}
It is closely related to the Lie algebra ${\mathfrak gl}_N$: the
exponents in the potential
correspond to the simple roots $e_i-e_{i+1}$, where $e_1,\ldots ,e_N$
denote the standard
basis elements in $\R^N$. More generally, the quantum Toda lattice
associated with a
real split semisimple (or reductive) Lie algebra $\mathfrak g$ with
Cartan subalgebra
$\mathfrak a$ has Hamiltonian given by
\[
\Delta_{\mathfrak a} - 2 \sum_{\alpha\in\Pi} d_\alpha e^{-\alpha
(x)} ,
\]
where $\Delta_{\mathfrak a}$ is the Laplacian on $\mathfrak a$,
$\Pi$ is a set of simple roots in ${\mathfrak a}^*$ and $d_\alpha$
are rational
numbers with a particular property \cite{glo}. For example, if
${\mathfrak g}=\mathfrak{so}_{2N+1}$, then we can identify $\mathfrak
a$ with $\R^N$,
take
\[
\Pi=\{e_1-e_2,e_2-e_3,\ldots ,e_{N-1}-e_N,e_N\},
\]
and the corresponding
Hamiltonian is given by
\[
\sum_{i=1}^N\frac{\partial^2}{\partial x_i^2} - 2 \sum_{i=1}^{N-1}
e^{x_{i+1}-x_i}
-e^{-x_N} .
\]
The connection between the (generalized) quantum Toda lattice and the
representation
theory of the corresponding Lie algebra $\mathfrak g$ was first
observed by Kostant \cite{k},
who showed that its eigenfunctions can be represented as particular
matrix elements
of infinite-dimensional representations of $\mathfrak g$. In the
simplest case when
${\mathfrak g}={\mathfrak sl}_2$ or ${\mathfrak gl}_2$, the
eigenfunctions are given in
terms of classical Whittaker functions (actually Macdonald functions).
For this reason,
they are often called $\mathfrak g$-Whittaker functions, or
$G$-Whittaker functions, if
${\mathfrak g}=\operatorname{Lie}(G)$. They also arise in the harmonic
analysis of
automorphic forms on Lie groups (see, e.g., \cite{bump}).
There is a spectral decomposition theorem in the general setting due to
Semenov--Tian--Shansky \cite{an}.
In this paper we will only consider the case ${\mathfrak g}={\mathfrak gl}_N$.
However, many of the constructions given throughout the paper have Lie-theoretic
interpretations and extend to the more general setting. This will be
indicated where appropriate.

The eigenfunctions of $H$ have the following integral
representation \cite{givental,jk,gklo,is}:
%
\begin{equation}\label{if}
\psi_\lambda (x) = \int_{\Gamma(x)} e^{\F_\lambda (T)} \prod
_{k=1}^{N-1}\prod_{i=1}^k dT_{k,i} ,
\end{equation}
where $\Gamma(x)$ denotes the set of real triangular arrays
$(T_{k,i}, 1\le i\le k\le N)$ with $T_{N,i}=x_i$, $1\le i\le N$,
and
\[
\F_\lambda (T)= \sum_{k=1}^N \lambda _k
 \Biggl( \sum_{i=1}^k T_{k,i}- \sum_{i=1}^{k-1} T_{k-1,i}  \Biggr)
-\sum_{k=1}^{N-1}\sum_{i=1}^k  ( e^{T_{k,i}-T_{k+1,i}} +
e^{T_{k+1,i+1}-T_{k,i}} ).
\]
This integral is defined for $\lambda\in\C^N$ and
has a recursive structure which we will now describe.
Write $H=H^{(N)}$, $\psi_\lambda =\psi^{(N)}_\lambda $. We will
drop these superscripts again later,
whenever they are unnecessary.
For convenience we define $H^{(1)}=d^2/dx^2$ and $\psi^{(1)}_\lambda
(x)=e^{\lambda x}$.
Following \cite{gklo}, for $N\ge2$ and $\theta\in\C$, define a
kernel on
$\R^N\times\R^{N-1}$ by
\[
Q^{(N)}_\theta(x,y) = \exp \Biggl( \theta
\Biggl ( \sum_{i=1}^N x_i- \sum_{i=1}^{N-1} y_i  \Biggr)
-\sum_{i=1}^{N-1}  ( e^{y_i-x_i} + e^{x_{i+1}-y_i} )
\Biggr) .
\]
Denote the corresponding integral operator by $\Q^{(N)}_\theta$,
defined on a suitable class of functions by
\[
\Q^{(N)}_\theta f(x)=\int_{\R^{N-1}} Q^{(N)}_\theta(x,y) f(y)\,dy.
\]
Then
%
\begin{equation}\label{ior}
\psi^{(N)}_{\lambda _1,\ldots ,\lambda _N} = \Q^{(N)}_{\lambda _N}
\psi^{(N-1)}_{\lambda _1,\ldots ,\lambda _{N-1}},
\end{equation}
and the integral formula (\ref{if}) can be re-expressed as
\[
\psi^{(N)}_\lambda =
\Q^{(N)}_{\lambda _N} \Q^{(N-1)}_{\lambda _{N-1}} \cdots\Q
^{(2)}_{\lambda _2} \psi^{(1)}_{\lambda _1} .
\]
Moreover, as remarked in \cite{gklo}, the following intertwining
relation holds:
%
\begin{equation}\label{i1}
\bigl(H^{(N)}-\theta^2\bigr)\circ\Q^{(N)}_\theta= \Q^{(N)}_\theta\circ H^{(N-1)}.
\end{equation}
This follows from the identity
\[
\bigl(H^{(N)}_x-\theta^2\bigr) Q^{(N)}_\theta(x,y) = H^{(N-1)}_y Q^{(N)}_\theta(x,y),
\]
which is readily verified. Combining (\ref{ior}) with the intertwining
relation (\ref{i1}) yields the
eigenvalue equation:
\[
H^{(N)} \psi^{(N)}_\lambda =  \Biggl( \sum_{i=1}^N \lambda
_i^2 \Biggr) \psi^{(N)}_\lambda .
\]

We note the following immediate consequences of the above integral
representation.
If $\lambda\in\iota\R^N$, then $\overline{\psi_\lambda(x)}=\psi
_{-\lambda}(x)$;
if $\lambda\in\iota\R^N$ and $\nu\in\R^N$, then $|\psi_{\lambda
+\nu}(x)|\le\psi_\nu(x)$.
It is also known (combining results from \cite{glo1,h,kl}) that for
each $x\in\R^N$,
$\psi_\lambda(x)$ is an entire function of $\lambda\in\C^N$.

The above construction has a representation-theoretic interpretation
which is
described in \cite{gklo}. It is closely related to the Gauss
decomposition and
has been extended to the other classical Lie algebras in \cite{glo}.
Encoded in the integrand are the
defining hyperplanes of the Gelfand--Tsetlin polytope associated with
the vector~$x$.

In the present setting (see, e.g., \cite{kl}), the spectral
decomposition theorem states
that the integral transform
%
\begin{equation}\label{u}
\hat f(\lambda)=\int_{\R^N} f(x)\psi_\lambda(x)\,dx
\end{equation}
defines an isometry from $L_2(\R^N,dx)$ into $L_2(\iota\R
^N,s_N(\lambda)\,d\lambda)$,
where $s_N(\lambda)\,d\lambda $ is the \textit{Sklyanin measure}
defined by
%
\begin{equation}\label{sklyanin}
s_N(\lambda )=\frac1{(2\pi\iota)^N N!} \prod_{j\ne k} \Gamma
(\lambda _j-\lambda _k)^{-1}.
\end{equation}

There is also a Mellin--Barnes type integral formula for $\psi_\lambda
$ due to Kharchev and
Lebedev \cite{kl,kl1,kl2} (see also \cite{is}). This is a kind of
dual of the above integral
representation and has a similar recursive structure. For $N\ge2$ and
$z\in\R$, define
a kernel on $\C^N\times\C^{N-1}$ by
\[
\hat{Q}^{(N)}_z(\lambda ,\gamma )= e^{z ( \sum\lambda _i-\sum
\gamma _i )} \prod_{i,j}\Gamma(\lambda _i-\gamma _j) .
\]
Then
%
\begin{equation}\label{kl}
\psi_\lambda ^{(N)} (x) = \int \hat{Q}^{(N)}_{x_1}(\lambda
,\gamma )
\psi_\gamma ^{(N-1)}(x_2,\ldots,x_N) s_{N-1}(\gamma )\,d\gamma
,\vspace*{-2pt}
\end{equation}
where the integral is along vertical lines with $\Re\gamma _i<\Re
\lambda _j$ for all $i,j$.
This construction also has a representation-theoretic interpretation
which is described
in \cite{gkl}. Gerasimov et al. \cite{glo1} give a clear account of
the nature of the duality
between the two constructions and, in particular, show how this duality yields
the following identity: for $\lambda,\nu\in\C^N$,
%
\begin{equation}\label{ri}
\int_{\R^N} e^{-e^{x_1-z}} \psi_\lambda (x) \psi_\nu(x)\,dx
=e^{z \sum(\lambda _i+ \nu_i)} \prod_{i,j}\Gamma(\lambda _i+\nu
_j) .\vspace*{-2pt}
\end{equation}
This is closely related to a Whittaker integral identity which was
conjectured by
Bump \cite{bump2} and later proved by Stade \cite{st2}, Theorem 1.1;
there is an
extensive literature on similar and related identities; see, for
example, \cite{bf,ja,st1}.

When $N=2$, the eigenfunctions $\psi_\lambda $ are given by
\[
\psi_\lambda (x) = 2 \exp \bigl(\tfrac12(\lambda _1+\lambda
_2)(x_1+x_2) \bigr)
K_{\lambda _1-\lambda _2} \bigl( 2 e^{(x_2-x_1)/2} \bigr),\vspace*{-2pt}
\]
where $K_\nu(z)$ is the Macdonald function. In this case,
the Givental's formula is equivalent to the integral formula
\[
K_\nu(z)=\frac12\int_0^\infty t^{\nu-1}\exp \biggl(-\frac
{z}2(t+1/t) \biggr)\,dt,\vspace*{-2pt}
\]
the contour integral representation (\ref{kl}) is equivalent to
\[
K_\nu(z)=\frac1{4\pi\iota}\int_{a-\iota\infty}^{a+\iota\infty
}\Gamma(s)\Gamma(s-\nu) \biggl(
\frac{z}2 \biggr)^{\nu-2s}\,ds,\qquad a>\max\{\Re{\nu},0\},\vspace*{-2pt}
\]
and the integral transform defined by (\ref{u}) is essentially (up to
a change
of variables) the Kontorovich--Lebedev transform.\vspace*{-2pt}

\section{The main result}

For $x,\nu\in\R^N$, denote by $\sigma^x_\nu$ the probability
measure on
the set $\Gamma$ of real triangular arrays $(T_{k,i})_{1\le i\le k\le N}$
defined by
\[
\int f \,d\sigma^x_\nu= \psi_\nu(x)^{-1}
\int_{\Gamma(x)} f(T)
e^{\F_\nu(T)} \prod_{k=1}^{N-1}\prod_{i=1}^k dT_{k,i}.\vspace*{-2pt}
\]
For $i=1,\ldots,N-1$, and continuous $\eta\dvtx (0,\infty)\to\R^N$, define
\[
(\T_i \eta) (t) = \eta(t) +  \biggl( \log\int_0^t e^{\eta
_{i+1}(s)-\eta_i(s)}\,ds \biggr) (e_i-e_{i+1}) ,\vspace*{-2pt}
\]
where $e_1,\ldots,e_N$ denote the standard basis vectors in $\R^N$.
Let $\Pi_1$ be the identity mapping ($\Pi_1\eta=\eta$)
and, for $2\le k\le N-1$, $\Pi_k = \T_1\circ\cdots\circ\T_{k-1}
\circ\Pi_{k-1}$.
Finally, we define
\[
\T= \Pi_N = (\T_1\circ\cdots\circ\T_{N-1})\circ\cdots\circ(\T
_1\circ\T_2) \circ\T_1 .\vspace*{-2pt}
\]
The main result of this paper is the following.\vadjust{\goodbreak}
%
\begin{theorem}\label{gmy}
{\smallskipamount=0pt
\begin{longlist}[(3)]
\item[(1)]
If $(W(t),t>0)$ is a standard Brownian motion in $\R^N$ with drift
$\nu$,
then $(\T W(t), t>0)$ is a diffusion process in $\R^N$ with
infinitesimal generator
given by
\[
\L_\nu= \frac12 \psi_\nu^{-1}  \Biggl(H-\sum_{i=1}^N\nu_i^2
\Biggr) \psi_\nu
= \frac12\Delta+ \nabla\log\psi_\nu\cdot\nabla.\vspace*{-2pt}
\]
\item[(2)]
For each $t>0$, the conditional law of $\{(\Pi_k W)_i(t), 1\le i\le
k\le N\}$,
given $\{\T W(s), s\le t; \T W(t)=x\}$, is given by $\sigma^x_\nu$.
\item[(3)] For each $t>0$, the conditional law of $W(t)$, given
$\{\T W(s), s\le t;\allowbreak \T W(t)=x\}$, is given by $\gamma^x_\nu$,
where
\[
\int_{\R^N} e^{(\lambda,y)} \gamma^x_\nu(dy) = \frac{\psi_{\nu
+\lambda}(x)}{\psi_\nu(x)},
\qquad\lambda\in\C^N.\vspace*{-2pt}
\]
\item[(4)]
If $\mu^{\nu}_t$ denotes the law of $\T W(t)$, then
\[
\mu^{\nu}_t(dx)=\exp \Biggl(-\frac12\sum_{i=1}^N\nu_i^2 t \Biggr)
\psi_\nu(x)\vartheta_t(x)\,dx,\vspace*{-2pt}
\]
where
%
\begin{equation}\label{theta}
\vartheta_t(x)=\int_{\iota\R^N} \psi_{-\lambda } (x) e^{\sum_i
\lambda _i^2 t/2} s_N(\lambda )\,d\lambda .\vspace*{-2pt}
\end{equation}
\end{longlist}}
\end{theorem}

It is easy to see that the process $(\{\Pi_k W)_i(t), 1\le i\le k\le
N\}, t>0)$ is Markov.
Indeed, setting $Z_{k,i}=(\Pi_k W)_i$, it follows from the
construction that $Z$ is a Markov
process taking values in $\Gamma$ which satisfies the system of stochastic
differential equations: $dZ_{1,1}=dW_1$ and, for $k=2,\ldots,N$,
\begin{eqnarray*}
dZ_{k,1} &=& dZ_{k-1,1} + e^{Z_{k,2}-Z_{k-1,1}}\,dt ,\\[-2pt]
dZ_{k,2} &=& dZ_{k-1,2} +  (
e^{Z_{k,3}-Z_{k-1,2}}-e^{Z_{k,2}-Z_{k-1,1}}  )\,dt, \\[-2pt]
& \vdots&\\[-2pt]
dZ_{k,k-1} &=& dZ_{k-1,k-1} +  (
e^{Z_{k,k}-Z_{k-1,k-1}}-e^{Z_{k,k-1}-Z_{k-1,k-2}}  )\,dt, \\[-2pt]
dZ_{k,k} &=& dW_k - e^{Z_{k,k}-Z_{k-1,k-1}}\,dt .\vspace*{-2pt}
\end{eqnarray*}
The infinitesimal generator of this process is given by
\[
\A_\nu=\frac12\sum_{1\le i\le k\le N} \frac{\partial^2}{\partial
z_{k,i}^2}
+\sum_{1\le i\le k<l\le N} \frac{\partial^2}{\partial
z_{k,i}\,\partial z_{l,i}}
+\sum_{1\le i\le k\le N} b_{k,i}(z)\, \frac\partial{\partial
z_{k,i}},\vspace*{-2pt}
\]
where
\begin{eqnarray*}
b_{1,1}(z)&=&\nu_1;
\\[-2pt]
b_{k,k}(z)&=&\nu_k-e^{z_{k,k}-z_{k-1,k-1}},\qquad k=2,\ldots,N;
\\[-2pt]
b_{k,1}(z)&=&e^{z_{k,1}-z_{k-1,1}},\qquad k=2,\ldots,N;
\\[-2pt]
b_{k,i}(z)&=&e^{z_{k,i+1}-z_{k-1,i}}-e^{z_{k,i}-z_{k-1,i-1}},\qquad
1<i<k\le N.
\end{eqnarray*}
The main content of Theorem~\ref{gmy} is the fact that $Z_{N,\cdot}$
is a Markov process,
with respect to its own filtration. The reason it holds is because
%
\begin{equation}\label{intertwining}
\L_\nu\circ\Sigma_\nu=\Sigma_\nu\circ\A_\nu,
\end{equation}
where $\Sigma_\nu$ is the Markov operator defined by
%
\begin{equation}\label{Sigma}
\Sigma_\nu f (x) = \int f(z) \sigma^x_\nu(dz).
\end{equation}
There is an additional (and nontrivial) technical issue related to the
fact that these processes start at a particular entrance law coming
from ``$-\infty$,'' but the
intertwining relation (\ref{intertwining}) lies at the heart of the
proof. Actually, the
proof of Theorem~\ref{gmy} given below is based on some intermediate
intertwining
relations which exploit the recursive structure of the quantum Toda
lattice and the
intertwining relation (\ref{intertwining}) is obtained as a
consequence, but it should,
nevertheless, be regarded as the analytic counterpart of Theorem~\ref{gmy}.
As far as we are aware, the intertwining relation (\ref{intertwining})
and its intermediaries
(given in Section~\ref{ir} below) have not been previously considered
in the literature.

The operator $\T$ was introduced (using a different notation) in the
paper \cite{oc03},
where it was surmised, based on heuristic arguments, that $\T W$ should
be a diffusion process
which has the same law as a Brownian motion conditioned, in an
appropriate sense,
on the asymptotic behavior of its exponential functionals. In \cite
{boc} it was observed
that such a conditioned Brownian motion can be defined and, moreover,
is closely related
to the quantum Toda lattice, thus providing the impetus for the present work.
The above notation used to define $\T$ follows a more general
framework which
has been developed in the papers \cite{bbo,bbo09}. It is shown
in \cite{bbo} that the
operators $\T_i$ satisfy the braid relations, that is,
\[
\T_i\circ\T_{i+1}\circ\T_i = \T_{i+1}\circ\T_i \circ\T
_{i+1},\qquad1\le i <N.
\]
It follows that for each element $\sigma\in{\mathfrak S}_N$ we can
uniquely define
\[
\T_\sigma= \T_{i_1}\circ\cdots\circ\T_{i_p},
\]
where $\sigma=(i_1,i_1+1)\cdots(i_p,i_p+1)$ is \textit{any} reduced
decomposition of $\sigma$
as a product of adjacent transpositions. The operator $\T$ corresponds
to the
longest element of ${\mathfrak S}_N$, that is, $\T=\T_{\sigma_0},$ where
\[
\sigma_0= \left(
\matrix{  1 & 2 & \cdots& N \cr N& N-1 & \cdots& 1
}
\right ) .
\]

The mapping
\[
\eta_{[0,t]} \mapsto\bigl(\{(\Pi_k \eta)_i(t), 1\le i\le k\le N\}, \{\T
\eta(s), s\le t\}\bigr)
\]
is a geometric variant of the RSK (Robinson--Schensted--Knuth) correspondence.
We will explain this connection later and give an interpretation of the measure
$\gamma_0^x$ appearing in Theorem~\ref{gmy} as a geometric analogue
of the
Duistermaat--Heckman measure associated with the point $x$.
The definition of the operator $\T$ extends naturally to other
Lie algebras, with ${\mathfrak S}_N$ replaced by the corresponding Weyl
group \cite{bbo,bbo09}.
It is natural to expect the analogue of Theorem~\ref{gmy} to hold in
this more general setting.

\section{The law of the partition function}

By Brownian scaling, it is easy to see that
the processes $(Z^N_t(\beta), t\ge0)$ and \hbox{$(\beta^{-2(N-1)}
Z^N_{\beta^2t}(1),t\ge0)$}
are identical in law, so for convenience we will define $Z^N_t=Z^N_t(1)$.
The transforma\-tion~$\T W$ is related to the random variable $Z^N_t$ as follows.
We first note that~$\T$ satisfies (cf. \cite{bbo}, Lemma 4.6)
%
\begin{equation}\label{sym}
(-\sigma_0)\circ\T=\T\circ(-\sigma_0),
\end{equation}
where $-\sigma_0 (\eta_1,\ldots,\eta_N)=(-\eta_N,\ldots,-\eta
_1)$ and $\eta_i$
denotes the $i${th} coordinate of the path $\eta$.
It is straightforward to see from the definition of $\T$ that
\begin{eqnarray*}
&&(\T W)_N(t)\\
&& \qquad =-\log\int_{0<s_1<\cdots<s_{N-1}<t}
e^{- (W_1(s_1)+W_2(s_1,s_2)+\cdots+ W_N(s_{N-1},t))}\,ds_1\cdots \,ds_{N-1}.
\end{eqnarray*}
From the relation (\ref{sym}) we have
\begin{eqnarray*}
&&(\T W)_1(t)\\
&&\qquad =\log\int_{0<s_1<\cdots<s_{N-1}<t}
e^{W_N(s_1)+W_{N-1}(s_1,s_2)+\cdots+ W_1(s_{N-1},t)}\,ds_1\cdots \,ds_{N-1}.
\end{eqnarray*}
Thus, if we set $W=(B_N,\ldots,B_1)$, then $\log Z^N_t = (\T W)_1(t)$
and we
deduce the following:
%
\begin{cor}\label{dp} The stochastic process $(\log Z^N_t,t>0)$ has
the same law as
the first coordinate of the diffusion process in $\R^N$ with
infinitesimal generator
\[
\L= \tfrac12 \psi_0^{-1} H \psi_0 = \tfrac12\Delta+ \nabla\log
\psi_0 \cdot\nabla,
\]
started according to the entrance law
\[
\mu_t(dx)=\psi_0(x)\vartheta_t(x)\,dx,\qquad t>0,
\]
where $\vartheta_t$ is given by (\ref{theta}). In particular, for
$u\in\R$, we have
\[
P(\log Z^N_t \le u) = \mu_t(\{x\in\R^N\dvtx  x_1 \le u\}) .
\]
\end{cor}

Note that the relation (\ref{sym}) also implies that the probability
measure $\mu_t$
is invariant under the transformation $-\sigma_0$. Combining Corollary
\ref{dp} with
the formula (\ref{ri}), we obtain [after shifting the contours in the
integral (\ref{theta})
to the left in order to apply Fubini's theorem] the following:
%
\begin{cor}\label{aaaaa} For $s>0$,
\[
Ee^{- s Z^N_t } =
\int s^{\sum\lambda _i} \prod_{i}\Gamma(-\lambda _i)^N e^{
(1/2)\sum_i\lambda _i^2 t} s_N(\lambda )\,d\lambda ,
\]
where the integral is along vertical lines with $\Re\lambda _i<0$ for
all $i$.
\end{cor}

The probability measure on $\iota\R^N$ with density proportional to
\[
e^{\sum_i \lambda _i^2 t/2} s_N(\lambda )\equiv\frac1{(2\pi\iota
)^N N!}
e^{\sum_i \lambda _i^2 t/2}\prod_{i>j}(\lambda _i-\lambda _j)\prod
_{i<j}\frac{\sin\pi(\lambda _i-\lambda _j)}{\pi}
\]
can be interpreted (up to a factor of $\iota\pi$)
as the law, at time $1/t$,
of the radial part of a Brownian motion in the symmetric space of
positive definite
Hermitian matrices or, equivalently, the law of the eigenvalues of a ``perturbed
GUE random matrix'' $A_N/\sqrt{t}+R_N/t$, where $A_N$ is an $N\times N$
GUE random matrix and $R_N$ is a diagonal matrix with entries given by the
vector $\pi(N-1,N-3,\ldots,1-N)$ (see, e.g., \cite{joc}).
In particular, it is a determinantal point process \cite{joh}. The above
expression for the moment generating function of $Z^N_t$ can thus be written
as a Fredholm determinant. It will be interesting to relate this,
in a suitable scaling limit, to the ``crossover distributions'' recently
introduced in the
context of KPZ and the stochastic heat equation by Sasamoto and
Spohn \cite{ss1,ss2,ss3,ss4}
and Amir, Corwin and Quastel \cite{acq}, building on recent work of
Tracy and
Widom \cite{tw1,tw2,tw3,tw4} on the asymmetric simple exclusion
process. See
also \cite{cdr,dot,dot2,dk,ps} for related recent developments.

We conclude this section by remarking that the other coordinates of~$\T
W(t)$ can
also be interpreted as logarithmic partition functions, as follows.
Define an ``up/right
path'' in $\R\times\Z$ to be an increasing path which either proceeds
to the right
or jumps up by one unit. For each sequence $0<s_1<\cdots<s_{N-1}<t$ we can
associate an up/right path $\phi$ from $(0,1)$ to $(t,N)$ which has
jumps between the
points $(s_i,i)$ and $(s_i,i+1)$, for $i=1,\ldots,N-1$, and is
continuous otherwise.
Then we can write
\[
(\T W)_1(t) = \log Z^N_t = \log\int e^{E(\phi)}\,d\phi,
\]
where
\[
E(\phi) = B_1(s_1)+B_2(s_2)-B_2(s_1)+\cdots+ B_N(t)-B_N(s_{N-1})
\]
and the integral is with respect to the Lebesgue measure on the
Euclidean set
of all such paths. There is an analogue of Greene's theorem in this
context~\cite{bbo}
(cf. \cite{kirillov})
which yields a similar formula for the other coordinates, namely, for
each $k=2,\ldots,N$,
\[
(\T W)_1(t)+\cdots+(\T W)_k(t) = \log\int e^{E(\phi_1)+\cdots
+E(\phi_k)}\,d\phi_1\cdots \,d\phi_k ,
\]
where the integral is with respect to the Lebesgue measure on the
Euclidean set of $k$-tuples\vadjust{\goodbreak}
of nonintersecting (disjoint) up/right paths with respective initial
points $(0,1),\ldots,(0,k)$ and
respective end points $(t,N-k+1)\allowbreak,\ldots,(t,N)$. An interesting
property of this formulation
is that it extends naturally to the continuum setting of KPZ and the
stochastic heat equation~\cite{OW}.

\section{The case $N=2$}

When $N=2$, the eigenfunctions $\psi_\nu$ are given by
\[
\psi_\nu(x) = 2 \exp \bigl(\tfrac12(\nu_1+\nu_2)(x_1+x_2) \bigr)
K_{\nu_1-\nu_2} \bigl( 2 e^{(x_2-x_1)/2} \bigr).
\]
In this case, Theorem~\ref{gmy} is equivalent to
the following theorem of Matsumoto and Yor \cite{my,my1}.
%
\begin{theo}
{\smallskipamount=0pt
\begin{longlist}[(3)]
\item[(1)]
Let $(B^{(\mu)}_t, t\ge0)$ be a standard one-dimensional Brownian
motion with drift~$\mu$,
and define
\[
Z^{(\mu)}_t=\int_0^t e^{2B^{(\mu)}_s-B^{(\mu)}_t}\,ds.
\]
Then $\log Z^{(\mu)}$ is a diffusion process with infinitesimal generator
\[
\frac12 \frac{d^2}{dx^2}+ \biggl(\frac{d}{dx}\log K_\mu
(e^{-x}) \biggr) \,\frac{d}{dx},
\]
where $K_\mu$ is the Macdonald function.
\item[(2)] The conditional law of $B^{(\mu)}_t$, given $\{Z^{(\mu)}_s,
s\le t; Z^{(\mu)}_t=z\}$,
is given by the generalized inverse Gaussian distribution
\[
\tfrac12 K_\mu(1/z)^{-1} e^{\mu x} \exp\bigl ( - \cosh(x)/z  \bigr)\,dx .
\]
\item[(3)]
The law of $Z^{(\mu)}_t$
is given by
\[
P\bigl(Z^{(\mu)}_t\in dz\bigr)= 2z^{-1}\theta_{1/z}(t)K_\mu(1/z)e^{-\mu^2
t/2}\,dz,
\]
where $\theta_r(t)$ is characterized by the Kontorovich--Lebedev transform
\[
2\int_0^\infty K_\lambda(r) \theta_r(t) \,\frac{dr}r = e^{\lambda
^2t/2},\qquad
\lambda\in\iota\R.
\]
\end{longlist}}
\end{theo}

The above Kontorovich--Lebedev transform can be inverted to obtain
\[
\theta_r(t)=\frac1{2\pi^2} \int_{-\iota\infty}^{\iota\infty}
K_\lambda (r) e^{\lambda ^2t/2} \lambda \sin(\pi\lambda )\,d\lambda
.
\]
The probability measure $H^{(1)}_r(dt)=I_0(r)^{-1}\theta_r(t)\,dt$ is
known as
the \textit{first  Hartman--Watson law} \cite{hw,my1}. It is also
characterized by
\[
\int_0^\infty e^{-\nu^2t/2} \theta_r(t)\,dt = I_\nu(r),\qquad\nu>0,
\]
where $I_\lambda$ is the modified Bessel function of the first kind.

\section{The zero-temperature limit}

By Brownian scaling, we can write down a version of Theorem~\ref{gmy}
for general $\beta>0$. We will state this in the case of zero drift.
For continuous $\eta\dvtx (0,\infty)\to\R^N$, define
%
\begin{eqnarray}
(\T^\beta_i \eta) (t) &=& \eta(t) + \frac1\beta\log \biggl( \beta
^2 \int_0^t
e^{\beta(\eta_{i+1}(s)-\eta_i(s))}\,ds \biggr) (e_i-e_{i+1}) ,\nonumber\\
\eqntext{i=1,\ldots,N-1;}
\\
\Pi^\beta_1&=&\mathit{Id}.;\qquad\Pi^\beta_k =
\T^\beta_1\circ\cdots\circ\T^\beta_{k-1} \circ\Pi^\beta
_{k-1}, \qquad2\le k\le N;
\nonumber\\
\T^\beta&=& \Pi^\beta_N = (\T^\beta_1\circ\cdots\circ\T^\beta
_{N-1})\circ\cdots
\circ(\T^\beta_1\circ\T^\beta_2) \circ\T^\beta_1 .\nonumber
\end{eqnarray}
Note that
\[
\frac1\beta\log Z^N_t(\beta) = (\T^\beta W)_1(t) - \frac
{N-1}{\beta}\log\beta^2.
\]
%
\begin{cor}\label{beta}
{\smallskipamount=0pt
\begin{longlist}[(3)]
\item[(1)]
If $W$ is a standard Brownian motion in $\R^N$, then $\T^\beta W$ is
a diffusion
in $\R^N$ with generator $\Delta/2
+ \nabla\log\psi_0(\beta \cdot ) \cdot\nabla$.
\item[(2)]
For each $t>0$, the conditional law of $\{(\Pi^\beta_k W)_i(t), 1\le
i\le k\le N\}$,
given $\{\T^\beta W(s), s\le t; \T^\beta W(t)=x\}$, is given by
$\sigma^{\beta x}_0(\beta \cdot )$.
\item[(3)] For each $t>0$, the conditional distribution of $W(t)$, given
$\{\T^\beta W(s),\allowbreak s\le t;  \T^\beta W(t)=x\}$, is given by $\gamma
^{\beta x}_0(\beta \cdot )$.
\item[(4)]
The law of $\T^\beta W(t)$ is given by $\mu_{\beta^2t}(\beta \cdot
 )$.
\end{longlist}}
\end{cor}

Letting $\beta\to\infty$, we recover the multidimensional version of
Pitman's ``$2M-X$'' theorem obtained in \cite{bj,oy,oc03a,bbo,bbo09}.
For continuous $\eta\dvtx (0,\infty)\to\R^N$, with $\eta(0)=0$, define
\begin{eqnarray*}
(\P_i \eta) (t) &=& \eta(t) -\inf_{0<s< t} \bigl(\eta_i(s)-\eta
_{i+1}(s)\bigr) (e_i-e_{i+1}) ,
\qquad i=1,\ldots,N-1;
\\
\Gamma_1&=&\mathit{Id}.;\qquad\Gamma_k =
\P_1\circ\cdots\circ\P_{k-1} \circ\Gamma_{k-1}, \qquad2\le k\le N;
\\
\P&=& \Gamma_N = (\P_1\circ\cdots\circ\P_{N-1})\circ\cdots\circ
(\P_1\circ\P_2) \circ\P_1 .
\end{eqnarray*}

By the method of Laplace, as $\beta\to\infty$, $\T^\beta W \to\P W$
uniformly on compact intervals and, for each $t>0$
and $1\le i\le k\le N$, $(\Pi^\beta_k W)_i(t) \to(\Gamma_k W)_i(t)$.
For $1\le k\le N$, $X^k=((\Gamma_k W)_1,\ldots,(\Gamma_k W)_k) $.
By construction, the stochastic process $\mathbb{X}(t)=(X^1(t),\ldots
,X^N(t))$, $t\ge0$,
is Markov and takes values in the Gelfand--Tsetlin cone
\[
\mathit{GT}_N=\{ (x^1,\ldots,x^N)\in\overline{C}_1 \times\cdots\times
\overline{C}_N\dvtx
x^{k+1}_{i+1} \le x_i^k \le x_i^{k+1},  1\le i\le k\le N-1\},
\]
where
\[
C_k=\{x\in\R^k\dvtx  x_1>\cdots>x_k\}.
\]
It is a $N(N-1)/2$-dimensional Brownian motion with singular covariance
reflected in
$\mathit{GT}_N$ via an explicit Skorohod reflection map. But we do not need
these facts, and
refer the reader to the papers \cite{oy,oc03a,bbo09} for details.

From the integral formula (\ref{if}) we have
\begin{eqnarray*}
&&\psi_0(\beta x)\\
&& \quad = \int_{\Gamma(\beta x)}
\exp\Biggl ( -\sum_{k=1}^{N-1}\sum_{i=1}^k  (
e^{T_{k,i}-T_{k+1,i}} + e^{T_{k+1,i+1}-T_{k,i}} )  \Biggr)
\prod_{k=1}^{N-1}\prod_{i=1}^k dT_{k,i}
\\
&& \quad = \beta^{N(N-1)/2}\\
&& \quad  \quad {}\times \int_{\Gamma(x)}
\exp \Biggl( -\sum_{k=1}^{N-1}\sum_{i=1}^k  \bigl( e^{\beta
(T'_{k,i}-T'_{k+1,i})} +
e^{\beta(T'_{k+1,i+1}-T'_{k,i})} \bigr)  \Biggr)
\prod_{k=1}^{N-1}\prod_{i=1}^k dT'_{k,i}.
\end{eqnarray*}
Write $x^k_i=T'_{k,i}$. As $\beta\to\infty$, if $x\in C_N$,
the integrand converges to 1 if $(x^1,\ldots,x^N)$ lies in the
Gelfand--Tsetlin polytope
\[
\mathit{GT}_N(x)=\{ (x^1,\ldots,x^N)\in \mathit{GT}_N\dvtx  x^N=x\},
\]
and 0 otherwise. It is well known (e.g., by Weyl's dimension
formula) that the
$N(N-1)/2$-dimensional Euclidean volume of $\mathit{GT}_N(x)$ is
\[
 \Biggl(\prod_{k=1}^{N-1} k!  \Biggr)^{-1} h(x),
\]
where $h$ is given by (\ref{h}).
It follows that, for $x\in C_N$,
%
\begin{equation}\label{asym0}
\lim_{\beta\to\infty} \beta^{-N(N-1)/2} \psi_0(\beta x)
=  \Biggl( \prod_{k=1}^{N-1} k!  \Biggr)^{-1} h(x) .
\end{equation}
Similarly, the probability measure
$\sigma^{\beta x}_0(\beta \cdot )$ converges as $\beta\to\infty$
to the uniform probability measure on $\mathit{GT}_N(x)$.
Putting all of this together, letting $\beta\to\infty$ in the
statement of
Corollary~\ref{beta}, we immediately recover parts (1) and (2) of the
following theorem.
%
\begin{theo}[(\cite{bj,oy,oc03a,bbo,bbo09})]\label{gp}
{\smallskipamount=0pt
\begin{longlist}[(3)]
\item[(1)]
If $W$ is a standard Brownian motion in $\R^N$, then $X^N=\P X$ is
a~Brownian motion conditioned (in the sense of Doob) never to exit
$C_N$.
\item[(2)]
The conditional law of $\mathbb{X}(t)$, given $\{X^N(s), s\le t;
X^N(t)=x\}$, is uniform on $\mathit{GT}_N(x)$.
\item[(3)]
The conditional law of $W(t)$, given $\{X^N(s), s\le t; X^N(t)=x\}$,
is given by the probability
measure $\kappa^x$ which is characterized by
\[
\int_{\R^N} e^{(\lambda,y)} \kappa^x(dy) =  \Biggl( \prod
_{k=1}^{N-1} k!  \Biggr)
\frac{\sum_{\sigma\in{\mathfrak S}_N} (-1)^\sigma e^{(\sigma
\lambda,x)} }{h(x)h(\lambda)}.
\]
\end{longlist}}
\end{theo}

Part (3) of the above theorem can be deduced from
part (2), noting that $\sum_{i=1}^k X^k_i = \sum_{i=1}^k W_i$, for each
$1\le k\le N$. Comparing this
with Corollary~\ref{beta}(3) yields the asymptotic formula: for $x\in C_N$,
%
\begin{equation}\label{asym1}
\lim_{\beta\to\infty} \beta^{-N(N-1)/2} \psi_{\lambda/\beta
}(\beta x)
= \frac{\sum_{\sigma\in{\mathfrak S}_N} (-1)^\sigma e^{(\sigma
\lambda,x)} }{h(\lambda)}.
\end{equation}
This formula can also be seen as a consequence of an alternative
representation of $\psi_\lambda$ as an alternating sum of \textit
{fundamental}
Whittaker functions \cite{h,kl,boc}.

The mapping
\[
\eta_{[0,t]} \mapsto\bigl(\{(\Gamma_k \eta)_i(t), 1\le i\le k\le N\}, \{
\P\eta(s), s\le t\}\bigr)
\]
is a continuous version of the RSK correspondence \cite{oc03a,bbo,bbo09}.
The mapping
\[
\eta_{[0,t]} \mapsto\bigl(\{(\Pi_k \eta)_i(t), 1\le i\le k\le N\}, \{\T
\eta(s), s\le t\}\bigr)
\]
is a continuous version of the geometric (or ``tropical'') RSK introduced by
Kirillov \cite{kirillov} (see also \cite{ny,b}).
The probability measure $\kappa^x$ is the (normalized)
Duistermaat--Heckman measure associated with the point $x$.
In this setting it can be interpreted, via the Harish--Chandra formula,
as the conditional distribution of the diagonal of a $N\times N$ GUE random
matrix given its eigenvalues $x$.
The probability measure $\gamma_0^x$ of
Theorem~\ref{gmy}
can thus be interpreted as a geometric analogue of the
Duistermaat--Heckman measure.
In keeping with this analogy, it is natural to record the following analogue
of the Littlewood--Richardson rule, which follows from Theorem~\ref{gmy}(3)
(cf.~\cite{bbo09}, Theorem 5.16(ii)).
For $s,t>0$, define $\tau_s W(\cdot)=W(s+\cdot)-W(s)$ and
\[
{\mathcal G}_{s,t}=\sigma\{\T W(r), 0< r\le s; (\T\tau_s W)(u),
0< u\le t\}.
\]

\begin{cor} For each $x,y\in\R^N$,
%
\begin{equation}\label{pf}
\frac{\psi_\lambda (x)}{\psi_0(x)}\frac{\psi_\lambda (y)}{\psi_0(y)}
= \int_{\R^N} \frac{\psi_\lambda (z)}{\psi_0(z)} \gamma^{x,y}(dz),
\end{equation}
where $\gamma^{x,y}$ is a probability measure on $\R^N$
which can be interpreted, for $s,t>0$, as the conditional law of $\T W(s+t)$
given $\mathcal G_{s,t}$, $\T W(s)=x$ and $(\T\tau_s W)(t)=y$.
\end{cor}

When $N=2$, (\ref{pf}) is equivalent to the formula
\[
K_\nu(z)K_\nu(w) = \frac12 \int_0^\infty e^{(-1/2) [t+(z^2+w^2)/t]}
K_\nu \biggl(\frac{zw}{t} \biggr)\,\frac{dt}{t}.
\]

Theorem~\ref{gp}, in the case $N=2$, is equivalent to Pitman's celebrated
``$2M-X$'' theorem \cite{pitman},
which states that, if $X_t$, $t\ge0$, is a standard one-dimensional
Brownian motion, then
$2\max_{0\le s\le t} X_s - X_t$, $t\ge0$, is a three-dimensional
Bessel process.
Setting $W=(B_N,\ldots,B_1)$ as before, the random variable\vadjust{\goodbreak}
$M^N_1$ defined by (\ref{mdef}) can be written as $M^N_1=X^N_1(1)$.
Thus, we also recover the fact \cite{bar,gtw} that $M^N_1$ has the
same law
as the largest eigenvalue of a $N\times N$ GUE random matrix.

Theorem~\ref{gp} has been generalized to arbitrary finite Coxeter
groups in the
papers \cite{bbo,bbo09}. The definition of the operator $\T$ also
extends naturally
to other Lie algebras, with ${\mathfrak S}_N$ replaced by the
corresponding Weyl group.
This is described in \cite{bbo,bbo09}, where various Lie-theoretic
interpretations are given.
It is natural to expect the analogue of Theorem~\ref{gmy} to hold in
this more general setting.

\section{Intertwining relations}\label{ir}

Consider the following extension of the operator $\Q^{(N)}_\theta$,
defined on a suitable class of functions $f\dvtx \R^N\times\R^{N-1}\to\R
$ by
\[
\qhat^{(N)}_\theta f(x)=\int_{\R^{N-1}} Q^{(N)}_\theta(x,y) f(x,y)\,dy .
\]
By a straightforward calculation, we obtain
%
\begin{equation}\label{i2}
\bigl(H^{(N)}-\theta^2\bigr)\circ\qhat^{(N)}_\theta= \qhat^{(N)}_\theta
\circ U^{(N)}_\theta,
\end{equation}
where
\begin{eqnarray*}
U^{(N)}_\theta &=& \sum_{i=1}^{N-1}\frac{\partial^2}{\partial y_i^2} -
2 \sum_{i=1}^{N-2} e^{y_{i+1}-y_i}
+ \sum_{i=1}^N\frac{\partial^2}{\partial x_i^2} \\
&&{}+ 2 ( \theta+e^{y_1-x_1}) \frac{\partial}{\partial x_1} \\
&&{}+ 2 ( \theta+e^{y_2-x_2}-e^{x_2-y_1} )  \frac{\partial}
{\partial x_2} \\
&&{}\;\; \vdots\\
&&{}+ 2 ( \theta+e^{y_{N-1}-x_{N-1}}-e^{x_{N-1}-y_{N-2}} )
\frac{\partial}{\partial x_{N-1}} \\
&&{} + 2( \theta - e^{x_N-y_{N-1}}) \frac{\partial}{\partial
x_N} .
\end{eqnarray*}
Further integration by parts yields
%
\begin{equation}\label{i3}
\bigl(H^{(N)}-\theta^2\bigr)\circ\qhat^{(N)}_\theta
= \qhat^{(N)}_\theta\circ V^{(N)}_\theta,
\end{equation}
where
\begin{eqnarray*}
V^{(N)}_\theta&=& \sum_{i=1}^{N-1}\frac{\partial^2}{\partial y_i^2}
- 2 \sum_{i=1}^{N-2} e^{y_{i+1}-y_i}
+ \sum_{i=1}^N\frac{\partial^2}{\partial x_i^2} \\
&&{}+ 2 \biggl ( \frac{\partial}{\partial y_1} + e^{x_2-y_1}  \biggr)\,
\frac{\partial}{\partial x_1} \\
&&{}+ 2  \biggl( \frac{\partial}{\partial y_2} + e^{x_3-y_2}-e^{x_2-y_1}
 \biggr)\, \frac{\partial}{\partial x_2} \\
&&{}\;\; \vdots\\
&&{}+ 2  \biggl( \frac{\partial}{\partial y_{N-1}} +
e^{x_N-y_{N-1}}-e^{x_{N-1}-y_{N-2}}  \biggr) \,\frac{\partial}{\partial
x_{N-1}} \\
&&{} + 2 ( \theta- e^{x_N-y_{N-1}} ) \,\frac{\partial
}{\partial x_N} .
\end{eqnarray*}
The intertwining relation (\ref{i3}) lies at the heart of this paper.

\section{\texorpdfstring{Proof of Theorem~\protect\ref{gmy}}{Proof of Theorem 3.1}}

We begin by using the intertwining relation~(\ref{i3}) to prove a
Markov functions result.
We will then proceed by induction to prove a version of Theorem~\ref
{gmy} for general
starting position. The final step will be to let the starting position
$x_0\to-\infty$ (in a sense
that will be made precise later).
Let $\nu\in\R^N$, and define
\[
\L^{(N)}_\nu= \frac12 \bigl(\psi^{(N)}_\nu\bigr)^{-1}  \Biggl( H^{(N)} - \sum
_{i=1}^N\nu_i^2  \Biggr) \psi^{(N)}_\nu.
\]
We consider a Markov process $((X(t),Y(t)),t\ge0)$ taking values in
$\R^N\times\R^{(N-1)}$,
defined as follows.
The process $Y$ evolves as an autonomous Markov process with
infinitesimal generator $\L^{(N-1)}_{\nu_1,\ldots,\nu_{N-1}}$. Let
$W$ be standard
one-dimensional Brownian motion, independent of $Y$, and define the
evolution of the
process $X$ via the stochastic differential equations
\begin{eqnarray*}
dX_1 &=& dY_1 + e^{X_2-Y_1}\,dt, \\
dX_2 &=& dY_2 +  ( e^{X_3-Y_2}-e^{X_2-Y_1}  )\,dt, \\
& \vdots&\\
dX_{N-1} &=& dY_{N-1} +  ( e^{X_N-Y_{N-1}}-e^{X_{N-1}-Y_{N-2}}
 )\,dt, \\
dX_N &=& dW +  ( \nu_N - e^{X_N-Y_{N-1}} )\,dt .
\end{eqnarray*}
Then $(X,Y)$ is a Markov process taking values in $\R^N\times\R^{(N-1)}$
with generator
\[
\G^{(N)}_\nu= \psi^{(N-1)}_{\nu_1,\ldots,\nu_{N-1}}(y)^{-1}
\Biggl ( V^{(N)}_{\nu_N} - \sum_{i=1}^{N-1} \nu_i^2  \Biggr) \psi
^{(N-1)}_{\nu_1,\ldots,\nu_{N-1}}(y) .
\]
Consider the Markov operator $\Lambda^{(N)}_\nu$ defined, for bounded
measurable functions
on $\R^N\times\R^{(N-1)}$, by
\[
\Lambda^{(N)}_\nu f (x) = \psi^{(N)}_\nu(x)^{-1} \int_{\R^{N-1}}
Q^{(N)}_\theta(x,y)
\psi^{(N-1)}_{\nu_1,\ldots,\nu_{N-1}}(y) f(x,y)\,dy .
\]
For $x\in\R^N$, define a probability measure $\lambda^x_\nu$ on $\R
^N\times\R^{(N-1)}$
by
\[
\int f \,d\lambda^x_\nu= \Lambda^{(N)}_\nu f(x).
\]
By (\ref{i3}), we have the intertwining relation
\[
\L^{(N)}_\nu\circ\Lambda^{(N)}_\nu= \Lambda^{(N)}_\nu\circ\G
^{(N)}_\nu.
\]
From the theory of Markov functions \cite{rp}, we conclude the following
proposition:
%
\begin{prop}\label{mf1}
Fix $x_0,\nu\in\R^N$ and let $(X,Y)$ be a Markov process with
infinitesimal generator $\G^{(N)}_\nu$,
started with initial law $\lambda^{x_0}_\nu$. Then $X$ is a~Markov
process with infinitesimal generator $\L^{(N)}_\nu$, started at
$x_0$. Moreover, for each $t\ge0$, the conditional
law of $Y(t)$, given $\{X(s), s\le t; X(t)=x\}$, is given by
\[
\psi^{(N)}_\nu(x)^{-1} Q^{(N)}_{\nu_N}(x,y) \psi^{(N-1)}_{\nu
_1,\ldots,\nu_{N-1}}(y)\,dy.
\]
\end{prop}

The next step is to deduce, by induction, an analogue of Theorem~\ref
{gmy} for general
starting position. We construct a Markov process $Z$ taking values in~%
$\Gamma$ as follows.
Let $W$ be a standard
Brownian motion in $\R^N$ with drift~$\nu$. The evolution of $Z$ is
defined recursively
by $dZ_{1,1}=dW_1$ and, for $k=2,\ldots,N$,
\begin{eqnarray*}
dZ_{k,1} &=& dZ_{k-1,1} + e^{Z_{k,2}-Z_{k-1,1}}\,dt, \\
dZ_{k,2} &=& dZ_{k-1,2} +  (
e^{Z_{k,3}-Z_{k-1,2}}-e^{Z_{k,2}-Z_{k-1,1}}  )\,dt, \\
& \vdots&\\
dZ_{k,k-1} &=& dZ_{k-1,k-1} +  (
e^{Z_{k,k}-Z_{k-1,k-1}}-e^{Z_{k,k-1}-Z_{k-1,k-2}}  )\,dt, \\
dZ_{k,k} &=& dW_k - e^{Z_{k,k}-Z_{k-1,k-1}}\,dt .
\end{eqnarray*}
%
\begin{prop}\label{mf2}
Fix $x_0,\nu\in\R^N$ and let $Z$ be the process defined as above
with initial law
$\sigma^{x_0}_\nu$.
Then $Z_{N,\cdot}$ is a Markov process with infinitesimal generator
$\L^{(N)}_\nu$, started at $x_0$.
Moreover, for each $t\ge0$, the conditional law of~$Z(t)$, given
$\{Z_{N,\cdot}(s), s\le t; Z_{N,\cdot}(t)=x\}$, is given by $\sigma
^x_\nu$, and
the intertwining relation (\ref{intertwining}) holds.
\end{prop}

Next we give a formula for the process $Z$ started at $Z(0)=z$ in terms
of the driving Brownian motion $W$.
For $i=1,\ldots,N-1$, and continuous~$\eta\dvtx\allowbreak(0,\infty)\to\R^N$, define
\[
(\T_i^\xi\eta) (t) = \eta(t) + \log \biggl( e^\xi+ \int_0^t
e^{\eta_{i+1}(s)-\eta_i(s)}\,ds  \biggr) (e_i-e_{i+1}) .\vadjust{\goodbreak}
\]
Fix $z\in\Gamma$ and, for $1\le i\le k\le N-1$, define
\[
\xi_{k,i}=z_{k,i}-z_{k+1,i+1}.
\]
Let $\Pi^z_1$ be the identity map and, for $2\le k\le N$,
\[
\Pi^z_k = (\T^{\xi_{k-1,1}}_1\circ\cdots\circ\T^{\xi
_{k-1,k-1}}_{k-1})\circ\Pi^z_{k-1} .
\]
Then, for $1\le i\le k\le N$, we can write
\[
Z_{k,i}(t) = z_{1,1} + (\Pi^z_k W)_i(t).
\]
For convenience we will write $\T^z=\Pi^z_{N}$ and note that
$Z_{N,\cdot}=z_{1,1}\mathbf{1}+\T^z W$, where $\mathbf
{1}=(1,1,\ldots,1)$.
Proposition~\ref{mf2} can now be restated as follows.
%
\begin{prop}\label{gsp}
Fix $x_0,\nu\in\R^N$. Let $W$ be a standard Brownian motion in $\R
^N$ with drift $\nu$
and $\zeta$ a random element of $\Gamma$ chosen according to the distribution
$\sigma^{x_0}_\nu$, independent of $W$. Then $Z_{N,\cdot}=\zeta
_{1,1}\mathbf{1}+\T^\zeta W$
is a Markov process with infinitesimal generator $\L^{(N)}_\nu$,
started at $x_0$.
Moreover, for each $t\ge0$, the conditional law of $Z(t)$, given
$\{Z_{N,\cdot}(s), s\le t; Z_{N,\cdot}(t)=x\}$, is given by $\sigma
^x_\nu$.
\end{prop}

For $k=1,\ldots,N$, define
\[
\rho^k= \biggl(\frac{k-1}2,\frac{k-1}2-1,\ldots,1-\frac
{k-1}2,-\frac{k-1}2 \biggr).
\]
We remark that the vector $\rho^k$ is half the sum of the positive
roots associated with the
Lie algebra $\mathfrak{gl}_k$.
To complete the proof of Theorem~\ref{gmy}, we will consider the
starting position
$x_0=-M\rho^N$, and let $M\to\infty$. For this we need to understand
the asymptotic
behavior of $\psi_\nu(-M\rho^N)$ and the probability measures
$\sigma^{-M\rho^N}_\nu$ as $M\to\infty$.
It was shown by Rietsch (\cite{rietsch}, Theorem 10.2) that the function
$-\F_0(T)$ on
$\Gamma(x)$ has a unique critical point $T^x$, which is a minimum, and
that the
Hessian is everywhere totally positive.
It is straightforward to verify from the critical point equations that
\[
\frac1k\sum_{i=1}^k T^x_{k,i} = \frac1N \sum_{i=1}^N x_i ,\qquad
1\le k\le N-1.
\]
Define ${\mathcal S}_\nu(T)=\F_\nu(T)-\F_0(T)$ and consider the
change of variables
\[
T'_{k,i}=T_{k,i}+M\rho^k_i, \qquad1\le1 \le k\le N.
\]
Then we can write
\[
\psi_\nu(-M\rho^N) = \int_{\Gamma(0)}
e^{{\mathcal S}_\nu(T')+e^{M/2}\F_0(T')} \prod_{k=1}^{N-1}\prod
_{i=1}^k dT'_{k,i}.
\]
It follows, by Laplace's method (see, e.g., \cite{es}, Theorem 4.14),
that the following asymptotic equivalence holds:
%
\begin{equation}\label{ae}
\psi_\nu(-M\rho^N) \sim C e^{-N(N-1)M/8} \exp\bigl(e^{M/2}\F_0(T^0)\bigr)
\end{equation}
as $M\to\infty$, where $C$ is a constant which is independent of $\nu$.
Moreover, recalling the above change of variables, we see that,
in probability, $\zeta_{k,i}-\zeta_{k+1,i+1}\to-\infty$ for each
$1\le i\le k\le N-1$ and $\zeta_{1,1}\to0$.
It follows by the continuous mapping theorem that
that $\zeta_{1,1}\mathbf{1}+\T^\zeta W$ converges in law to~$\T W$,
and, for each $t>0$, $\{(\Pi^\zeta_k W)_i(t), 1\le i\le k\le N\}$
converges in
law to $\{(\Pi_k W)_i(t), 1\le i\le k\le N\}$.
We conclude that $\T W$ is a diffusion with generator $\L^{(N)}_\nu$,
and that the conditional law of $\{(\Pi_k W)_i(t), 1\le i\le k\le N\}$,
given $\{\T W(s), s\le t; \T W(t)=x\}$, is $\sigma^x_\nu$.
This proves parts (1) and (2) of the theorem.
Part (3) of the theorem follows from part~(2), noting that for each $k\le N$,
\[
W_k = \sum_{i=1}^k (\Pi_k W)_i - \sum_{i=1}^{k-1} (\Pi_{k-1} W)_i .
\]
Part (4) follows from part (3) by the spectral decomposition theorem.
%
\begin{rem}
The asymptotic equivalence (\ref{ae}) is well known in the case $N=2$ and
can be compared to the full asymptotic expansion obtained in~\cite{bh}
in the
case $N=3$, where it was remarked that the leading term in the
expansion is
independent of the parameter $\nu$.
\end{rem}

\section{\texorpdfstring{A symmetric version of Proposition~\protect\ref{mf2}}{A symmetric version of Proposition 8.2}}

Proposition~\ref{mf2} has a ``symmetric'' analogue which can be regarded
as a
geometric version of a result of Dubedat \cite{dub} in the case $N=2$,
and Warren \cite{w} in the general case. It is obtained by applying the
intertwining relation (\ref{i2}) rather than (\ref{i3}). In this
case, we construct a
Markov process $S$ on $\Gamma$ as follows.
Let $\{W_{k,i}, 1\le i\le k\le N\}$ be a collection of independent standard
one-dimensional Brownian motions. The evolution of $S$ is
defined recursively by $dS_{1,1}=dW_{1,1}$ and, for $k=2,\ldots,N$,
\begin{eqnarray*}
dS_{k,1} &=& dW_{k,1} + (\nu_k+e^{S_{k-1,1}-S_{k,1}})\,dt, \\
dS_{k,2} &=& dW_{k,2} + ( \nu_k+e^{S_{k-1,2}-S_{k,2}}-e^{S_{k,2}-
S_{k-1,1}} )\,dt, \\
& \vdots& \\
dS_{k,k-1} &=& dW_{k,k-1} + ( \nu_k+ e^{S_{k-1,k-1}-S_{k,k-1}}-e^
{S_{k,k-1}-S_{k-1,k-2}} )\,dt, \\
dS_{k,k} &=& dW_{k,k} + ( \nu_k - e^{S_{k,k}-S_{k-1,k-1}} )\,dt .
\end{eqnarray*}
%
\begin{prop}\label{mf3}
Fix $x_0,\nu\in\R^N$ and let $S$ be the process defined as above
with initial law
$\sigma^{x_0}_\nu$.
Then $S_{N,\cdot}$ is a Markov process with infinitesimal generator
$\L^{(N)}_\nu$, started at $x_0$.
Moreover, for each $t\ge0$, the conditional law of $S(t)$, given
$\{S_{N,\cdot}(s), s\le t; S_{N,\cdot}(t)=x\}$, is $\sigma^x_\nu$.
\end{prop}

In the case $N=2$, with zero drift, we deduce the following corollary:
%
\begin{cor}
Let $B_1,B_2$ and $B_3$ be independent standard one-dimen\-sional
Brownian motions.
Define
\begin{eqnarray*}
X(t)&=& B_1(t) + \log\int_0^t e^{B_2(s)-B_1(s)}\,ds,
\\
Y(t)&=& B_3(t) - \log\int_0^t e^{B_3(s)-B_2(s)}\,ds.
\end{eqnarray*}
Then $(X+Y)/\sqrt{2}$ is a standard Brownian motion and
$(X-Y)/\sqrt{2}$ is a~diffusion process [independent of
$(X+Y)/\sqrt{2}$]
with infinitesimal generator
\[
\frac12 \frac{d^2}{dx^2}+ \biggl(\frac{d}{dx}\log K_0(e^{-x}) \biggr)
\,\frac{d}{dx}.
\]
\end{cor}

\section*{Acknowledgments} The author would like to thank the referees
for careful
reading and helpful suggestions which have led to an improved version
of the paper.
Thanks also to Philippe Biane,
Alexei Borodin,
 Philippe Bougerol, Ivan Corwin, Timo
Sepp\"{a}l\"{a}inen
and Herbert Spohn for helpful and stimulating discussions, and Dan Bump for
pointers to the literature on Whittaker integral identities.


%

\printaddresses

\end{document}